\documentclass[12pt,a4paper]{article}
\usepackage[T2A]{fontenc}
\usepackage[cp1251]{inputenc}
\usepackage{amsmath,amssymb,amsthm,amsfonts,amscd}
\begin{document}
\newtheorem{example}{Example}
\newtheorem{proposition}{Proposition}
\newtheorem{lemma}{Lemma}
\newtheorem{remark}{Remark}
\newtheorem{theorem}{Theorem}
\newtheorem{definition}{Definition}
\def\K{{\bf K}}
\def\SH{{S\Bbb{H}}}
\def\SO{{S\Bbb{O}}}
\def\SU{{S\Bbb{U}}}
\def\S{{\Bbb{S}}}
\def\O{{\Bbb{O}}}
\def\R{{\Bbb R}}
\def\I{{\Bbb I}}

\def\Z{{\Bbb Z}}
\def\P{{\Bbb P}}
\def\RP{{\Bbb R}\!{\rm P}}
\def\N{{\Bbb N}}
\def\C{{\Bbb C}}
\def\H{{\Bbb H}}
\def\Q{{\bf Q}}
\def\A{{\bf A}}
\def\D{{\bf D}}
\def\k{{\bf k}}
\def\K{{\bf K}}
\def\B{{\bf B}}
\def\E{{\bf E}}
\def\F{{\bf F}}
\def\V{\vec{\bf V}}
\def\L{{\bf L}}
\def\G{{\bf G}}
\def\c{{\bf c}}
\def\i{{\bf i}}
\def\j{{\bf j}}

\def\fr{{\operatorname{fr}}}
\def\st{{\operatorname{st}}}
\def\mod{{\operatorname{mod}\,}}
\def\cyl{{\operatorname{cyl}}}
\def\dist{{\operatorname{dist}}}
\def\grad{{\bf{grad}}}
\def\div{{\operatorname{div}}}
\def\rot{{\operatorname{rot}}}

\def\R{{\Bbb R}}
\def\B{{\bf B}}
\def\e{{\bf e}}
\def\L{{\bf L}}
\def\valpha{\vec{\alpha}}
\def\vxi{\vec{\xi}}
\sloppy
\title{Higher Helicity of Magnetic Lines and Arf-invariants}
\author{Petr M. Akhmetiev} 
\date{IZMIRAN,  Troitsk, Moscow}
\maketitle
\medskip
\bigskip
\sloppy
\section*{Introduction}

V.I.Arnold formulated the following problem [\cite{Arn}, Problem 1984-
12]: "To transform the asymptotic ergodic definition of the Hopf invariant
of divergence-free vector fields to the theory of S.P.Novikov which generalizes the Whitehead product of homotopy groups of spheres"'. In the
paper we recall and simplify (a partial) solution of the problem from the \cite{A4} and present new results, which generalize the problem to non-simply connected manifolds.

In the first section, we present an additional motivation of the Arnold Problem, which is based on mean magnetic field theory.
We use geometrical considerations due to K.Moffatt and formulate properties of invariants in ideal MHD, which are
asymptotic and ergodic properties.  

Then we recall the definition of the quadratic helicity invariant and of the higher asymptotic ergodic $M$-invariant.
We present a  simpler new proof (in part) that the $M$-invariant is ergodic. The $M$-invariant  is a higher invariant, this means that for the magnetic field with closed magnetic lines the invariant is not a function of 
pairwise linking numbers of the magnetic lines. This property is based of the following fact: an arithmetic residue of the $M$-invariant for a triple of closed magnetic lines, which is a model of a  link with 
even pairwise linking numbers,
 coincides with the Arf-invariant (about the Arf-invariant, or, the Rokhlin-Robertello invariant, see \cite{G-M}).

The new results concern  magnetic fields on closed 3-dimensional manifolds and use the $M$-invariant.
The manifolds with magnetic field, that we consider   are not, generally speaking, simply-connected. This manifold is
assumed homogeneous and  is a rational Poncar\'e sphere.
One can try to transform results on the asymptotics and ergodicity of the $M$-invariant  for the magnetic fields on the standard sphere $S^3$ to an arbitrary rational homology sphere $\Sigma$. To make  this idea precise we generalize
the Arf-invariants of classical semi-boundary links  (including the Arf-Brown $\Z/8$-invariant) (see \cite{G-M}) and we introduce a new
Arf-invariant, called the hyperquaternionic Arf-invariant.      

This generalization could clarify the relationship between the $M$-invariant and homotopy groups of spheres. 
It is well-known that the helicity invariant is a specification of the Hopf invariant, see \cite{A-Kh} for details.
The Hopf invariant determines the homotopy group $\pi_3(S^2)$, the stabilization of this homotopy group is denoted by $\Pi_1$. The group $\Pi_1$ contains the only non-trivial element with the Hopf invariant one. 

The Arf-invariant describes the stable homotopy group $\Pi_2$ via the  geometrical approach due to L.S.Pontrjagin.
The Arf-Brown invariant describes the 2-torsion of the stable homotopy group 
$\Pi_3$, this result follows from  V.A.Rokhlin's theorems.  The hyperquaternionic Arf-invariant 
describes the 2-torsion of the stable homotopy group $\Pi_7$. This group was calculated by J.P.Serre using the algebraic  approach. The complexity of the $M$-invariant relates with the fundamental group of rational homology sphere $\Sigma$.

The results were presented at the A.B.Sossinsky Topological Seminar in IMU September-October 2004. 
A preliminary result was presented at the conference on differential equations,  organized by V.P.Leksin in Kolomna, June 2014.
The author was supported in part by RFBR grants 15-02-01407, 15-01-06302.


\section*{The mean magnetic field equation}

Let us consider, as in \cite{R}, the domain  
$\Omega$ in $\R^3$, which is compact for simplicity, with a conductive liquid. In $\Omega$ a velocity field $\bf{u}$
of the liquid and a magnetic field $\bf{B}$ are well-defined. Moreover, the following decomposition
of the considered vector-fields into a mean part and a random part is well defined:
$$ \bf{B} = \bar{\bf{B}} + \bf{B'}; \qquad \bf{u} = \bar{\bf{u}} + \bf{u'}. $$

Assume that the mean velocity field 
 $\bar{\bf{u}}(t)$ is done, then the equation for the mean magnetic field is following: 
\begin{eqnarray}\label{dynamo}
\begin{array}{c}
\rot (\eta \rot \bar{\bf{B}}) - \rot( \bar{\bf{u}} \times \bar{\bf{B}} + \E) + \frac{\partial{\bar{\bf{B}}}}{\partial t} = 0, \\
\E = \overline{\bf{B'}  \times \bf{u'}},  \quad \div(\bar{\bf{B}})=0. \\
\end{array}
\end{eqnarray}
The equation
 (\ref{dynamo}) is called the kinematic dynamo equation. Assuming  $\eta=0, \E=0$ 
this equation means that the magnetic field is frozen-in. 

Assume that the following equation is satisfied:
\begin{eqnarray}\label{eds}
 \E  = \alpha \bar{\bf{B}} - \beta \rot(\bar{\bf{B}}). 
\end{eqnarray} 
Then, using the condition that $\alpha$ changes the sign with respect to the mirror symmetry and using additional
simplifing assumptions we get:  
\begin{eqnarray}\label{alpha}
\alpha \sim   \overline{(\bf{u'},\rot(\bf{u'}))}, 
\end{eqnarray}
where the function
$(\bf{u'},\rot(\bf{u'}))$ is called the density of (a small-scaled) the hydrodynamic helicity.
Denote the hydrodynamic helicity by 
 $\chi_{\bf{u}'} = \int (\bf{u}',\rot(\bf{u}')) d\Omega$.

Take the scalar product of the both sides of the equation
 (\ref{eds}) with the vector  $\bar{\bf{B}}$, assuming for simplicity that  $\eta=0$, and take the integral over the domain $\Omega$. We get, using  $\E = \frac{\partial \A}{\partial t}$, the equation, which describes the transport of the magnetic helicity  $\chi_{\bar{\B}}= \int (\A,\bar{\B}) d\Omega$:
\begin{eqnarray}\label{helicity}
\frac{d\chi_{\bar{\B}}}{dt} = 2\alpha \int(\bar{\B},\bar{\B})d\Omega - 2\beta \int (\bar{\B},\rot(\bar{\B})) d\Omega.
\end{eqnarray}
The integral $U_{\bar{\B}}=2\int(\bar{\B},\bar{\B})d\Omega$ is called the magnetic energy (of the mean field), the integral 
$\chi_{\rot \bar{\B}}=2\int (\bar{\B},\rot(\bar{\B})) d\Omega$ 
is called the current helicity (of the mean field).

\subsection*{Topological considerations concerning the transport equation of the magnetic helicity}

In the paper \cite{M} by K.Moffatt the equation (\ref{helicity}) is discussed from point of view of geometry 
of magnetic lines. Assume that a support of a magnetic field consists of a finite set of magnetic tubes, see \cite{B-F}. This means that the magnetic fields $\bar{\bf{B}}$,  $\bf{B'}$  are inside the tubes and is tangent to the surfaces of the tubes.
Additionally, assume that the same collection of the tubes is a support of a velocity field $\bf{u}'$.
With this assumption the vorticity  field points along the central axis of the each tube.

The magnetic and hydrodynamic tubes one may define such that the following condition, which is called "`force-free"' condition is satisfied: 
 $\rot \bf{u}' \sim \bf{u}'$, $\rot \B' \sim \B'$.

With the considered assumption it is not hard to prove, using the formula $(\ref{alpha})$,  that  the mean magnetic field $\bar{\bf{B}}$ in the collection of tubes tends $\vert \alpha \vert $-exponentially
to $+\infty$, if the absolute value of $\alpha$ is sufficiently large. This is called the $\alpha$-effect.

The contribution of the second term in the right side of the equation  
 (\ref{helicity}) is given by the Calugareanu formula, see \cite{M-R}.
The magnetic helicity inside the only magnetic tube is calculated by the formula: 
 $$ \chi_{\B} = \Phi^2 Lk, $$
where $\Phi$ is the integral magnetic flow trough a transversal section of the tube, $Lk$ is the self-linking number of the magnetic tube. For two magnetic tubes  $\Omega_1$, $\Omega_2$ with magnetic flows $\Phi_1$, $\Phi_2$ the magnetic helicity is calculated by the formula: 
$$ \chi_{\B} = \Phi_1^2 Lk(1,1) + 2 \Phi_1 \Phi_2 Lk(1,2) +  \Phi_2^2 Lk(2,2). $$

\subsection*{Magnetic tubes}

 \vspace{10mm}
\protect
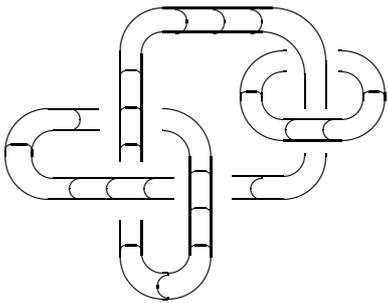
\begin{picture}(350,80)

\put(95,42){\oval(34,34)[l]} \put(95,42){\oval(18,18)[l]}

\put(115,42){\oval(34,34)[r]} \put(115,42){\oval(18,18)[r]}

\put(25,20){\oval(70,34)[l]} \put(50,-5){\oval(34,60)[b]}
\put(50,45){\oval(34,60)[tl]} \put(75,20){\oval(70,34)[br]}

\put(25,20){\oval(50,18)[l]} \put(50,-5){\oval(18,40)[b]}
\put(50,45){\oval(18,42)[tl]} \put(75,20){\oval(54,17)[br]}

\put(50,37){\oval(120,76)[tr]}

\put(95,33){\line(1,0){20}} \put(95,25){\line(1,0){20}}
\put(102,37){\line(0,1){12}}

\put(62,49){\oval(80,34)[tr]}

\put(49,-5){\oval(36,84)[tr]} \put(49,-5){\oval(20,68)[tr]}

\put(33,17){\line(0,1){28}}

\put(41,17){\line(0,1){28}}

\put(25,11){\line(1,0){28}} \put(25,3){\line(1,0){28}}

\put(50,66){\line(1,0){20}}

\put(-5,20){\oval(10,7)[t]} \put(37,20){\oval(8,7)[t]}

\put(82,40){\oval(8,7)[t]} \put(128,40){\oval(8,7)[t]}

\put(99,29){\oval(8,8)[l]} \put(113,29){\oval(8,8)[l]}

\put(37,34){\oval(8,7)[t]} \put(37,48){\oval(8,7)[t]}
\put(63,10){\oval(8,7)[t]} \put(63,-4){\oval(8,7)[t]}
\put(63,-18){\oval(8,7)[t]} \put(37,-18){\oval(8,7)[t]}

\put(14,33){\oval(8,8)[r]} \put(54,70){\oval(8,9)[r]}
\put(68,70){\oval(8,9)[r]} \put(82,70){\oval(8,9)[r]}

\put(18,7){\oval(8,8)[l]} \put(32,7){\oval(8,8)[l]}
\put(46,7){\oval(8,8)[l]} \put(86,7){\oval(8,8)[l]}
\put(51,-30){\oval(8,9)[l]}

 \protect
\end{picture}

\vspace{35mm}


 \vspace{10mm}
\protect
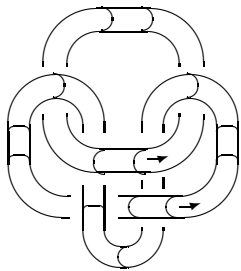
\begin{picture}(350,80)


\put(210,53){\oval(60,44)[t]} \put(185,28){\oval(36,44)[t]}
\put(235,28){\oval(36,44)[t]} \put(210,35){\oval(60,44)[b]}
\put(210,-8){\oval(30,30)[b]} \put(190,28){\oval(46,64)[bl]}
\put(230,28){\oval(46,64)[br]}

\put(210,53){\oval(42,26)[t]} \put(185,28){\oval(20,26)[t]}
\put(235,28){\oval(20,26)[t]} \put(210,35){\oval(42,26)[b]}
\put(210,-8){\oval(14,14)[b]} \put(190,28){\oval(30,48)[bl]}
\put(230,28){\oval(30,48)[br]}

\put(195,-8){\line(0,1){16}}

\put(203,-8){\line(0,1){16}}

\put(208,4){\line(1,0){24}}

\put(208,-4){\line(1,0){25}}

\put(225,7){\line(0,1){4}}

\put(217,7){\line(0,1){3}}

\put(158,30){\makebox(0,0)}
\put(246,74){\makebox(0,0)}

\put(171,27){\oval(8,7)[t]} \put(171,16){\oval(8,7)[t]}
\put(249,27){\oval(8,7)[t]} \put(249,16){\oval(8,7)[t]}
\put(199,-3){\oval(8,7)[t]}

\put(202,71){\oval(8,9)[r]} \put(217,71){\oval(8,9)[r]}
\put(184,46){\oval(8,9)[r]} \put(234,46){\oval(8,9)[r]}

\put(216,0){\oval(8,8)[l]} \put(230,0){\oval(8,8)[l]}
\put(203,17){\oval(8,9)[l]} \put(218,17){\oval(8,9)[l]}
\put(212,-19){\oval(8,8)[l]}

\put(235,0){\vector(3,1){4}} \put(223,18){\vector(3,1){4}}

\put(231,0){\line(1,0){4}} \put(219,18){\line(1,0){4}}

 \protect
\end{picture}

\vspace{35mm}

Assume that magnetic lines in each magnetic tube on the figure are untwisted with respect to the plane of the projection. Then for the first pair of magnetic tubes  we have:
$Lk(1,1)=-3$, $Lk(1,2)=\pm1$, $Lk(2,2)=0$; for the second pair of magnetic tubes we have:   $Lk(1,1)=-2$, $Lk(1,2)=Lk(2,2)=0$.

Assume that a magnetic field is inside the only magnetic tube. In this case the self-helicity of the magnetic tube $\Omega$
is calculated by the Calugareanu formula:
 $$Lk = Wr + Tw. $$
(for the simplicity we assume that the magnetic flow trough the magnetic tube is equal to $1$).
Consider the central line $l$ of the magnetic tube 
 $\Omega$ and consider the vector field $\xi$, which is perpendicular to $l$, 
and points from $l$ to anyone magnetic line in $U$. Let us say that the magnetic tube $\Omega$ is "`hairless"'
if the derivative of the field along the tangent vector of $l$ is trivial at each point of $l$. 

For an arbitrary magnetic tube $\Omega$ with a central line $l$ consider the "`hairless"' magnetic tube $\Omega_0$ with the same central line and the same (unite) magnetic flow. The self-linking number for $\Omega_0$ is calculated by the 
Calugareanu formula as the sum of the self-linking number $Wr$ of $\Omega_0$ with the twisting number $Tw$ of the magnetic tube $\Omega$ with respect to $\Omega_0$.

It is not hard to prove the following equation:
$$ \chi_{\rot \B} = C \Phi^2 Tw, $$ 
where $\chi_{\rot \B}$ is the current helicity of the magnetic field in $\Omega$, $Tw$ is the twisting number of the magnetic tube $\Omega$, $C$ is the coefficient of an order $1$, which depends of a geometry of magnetic lines inside $U$.  
Therefore the second term in the formula 
(\ref{helicity}) detects an untwist of the magnetic tube, see \cite{A-K-K} for more details. 

\section*{Quadratic helicity and ergodic integrals}

Assume for a simplicity that magnetic lines in a magnetic tube $\Omega$ are closed.
The magnetic tube  $\Omega$ is characterized by the magnetic helicity integral, this integral is equal to the mean pairwise linking number of  magnetic lines is the magnetic tube $\Omega$, which is normalized by magnetic flows thought the collection of infinitesimal magnetic lines $\Omega$.

The magnetic tube $\Omega$ is also characterized by various combinatorial invariants 
$I(L_1,L_2,L_3)$, which are calculated for various collections   $\{L_1,L_2,L_3\}$ of $k$ magnetic lines (we assume $k=3$
for simplicity). In this case we may assume that lines of  collections  are inside the magnetic tubes $\Omega_1, \Omega_2, \Omega_3$ correspondingly, some magnetic tubes could coincide. In a particular interesting case we have the only magnetic tube, magnetic lines of the collections are inside of this tube. 

What are required conditions for a combinatorial invariant $I$, which can be apply to describe  magnetic fields?  
From the consideration above of the equation (\ref{helicity}) we have to assume the following conditions:

$\bullet$ C1. The invariant  $I$ is of a finite-type invariant of an order $t$ in the sense of V.A.Vassiliev.

$\bullet$ C2. The invariant  $I$ is characterized by a positive integer $s$, which is called the asymptotic denominator. Take a link $(L_1,L_2,L_3)$, which is formed by central lines of disjoint magnetic tubes $\Omega_1,\Omega_2, \Omega_3$. Denote by $(rL_1,rL_2,rL_3)$  the $r$-time spinning link, which is constructed from $(L_1,L_2,L_3)$ by the $r$-fold 
spinning along the central line of the  corresponding magnetic tubes $\Omega_1, \Omega_2, \Omega_3$. The following equation is satisfied:
$$r^{3s}I(L_1,L_2,L_3)  = I(rL_1,rL_2,rL_3) + O(r^{3s-1}). $$

$\bullet$ C3. Assume we have two disjoint magnetic tubes  $\Omega_2, \Omega_3$ and we have two parallel magnetic lines,
which is a 2-component link $(L_1, L_2)$ in $\Omega_2$, and a magnetic line $L_3$ is a central line in $\Omega_3$. Take a magnetic tube $\Omega_2^{tw}$, which is obtained from the magnetic tube $\Omega_2$ by a twist, $\Omega_2 \mapsto \Omega_2^{tw}$, $Tw(\Omega_2) = Tw(\Omega_2^{tw}) + const$. Take two parallel magnetic lines
$L_1^{tw}, L_2^{tw}$ in $\Omega_2^{tw}$. Take the $r$-time spinning link $(rL_1,rL_2)$, each commponent of this link  is rotated along the central line $L_1 =L_2$ of $\Omega_2$ $r$ times. Take the $r$-time spinning link $(rL^{tw}_1,rL^{tw}_2)$, each commponent of this link  is rotated along the central line $L^{tw}_1 =L^{tw}_2$ of $\Omega^{tw}_2$ $r$ times. Take 3-component links $(rL_1,rL_2,rL_3)$, $(rL^{tw}_1,rL^{tw}_2,rL_3)$. The following formula is satisfied:
$$I(rL_1,rL_2,rL_3)  - I(rL^{tw}_1,rL^{tw}_2,rL_3) = O(r^{3s}). $$

$\bullet$ C4 (Condition-Definition). Assume that the invariant $I$ is not a function of 
pairwise linking numbers of components of the link (for a 3-component link we get 3 pairwise linking numbers).
In this case we say that $I$ is a higher invariant. 

\subsubsection*{Quadratic helicity}

We shell give an non-formal definition. A precise definition of the quadratic helicity (without the assumption that a magnetic line is  non-closed) is presented in \cite{A}.
Consider a finite collection 
$\{l_i\}$ of  $N$  magnetic lines (each line is closed for simplicity).  $1 \le i \le N$, $N >> 1$, and consider a symmetric $N \times N$ matrix 
with zero elements on the main diagonal, which is defined  by the pairwise linking numbers   $lk(i,j), \quad 1 \le i < j \le N$ of magnetic lines.

Define the quadratic helicity integral 
 $\chi^{(2)}$ by the formula: 
\begin{eqnarray}\label{11}
\chi^{(2)} = \Phi_i\Phi_j\Phi_k\sum_{i,j,k} lk(i,j)lk(j,k), \quad i \ne j, j \ne k, k \ne i, \quad 1<i,j,k<N, 
\end{eqnarray}
where $\Phi_i, \Phi_j, \Phi_k$ are integral magnetic flows trough infinitesimal thin magnetic tubes $\Omega_i$, $\Omega_j$, $\Omega_k$, which are formed by the corresponding magnetic lines.

The quadratic helicity satisfy Conditions C1,C2,C3 for 
$s=\frac{4}{3}$, $t=2$. It is interesting to remark that 
\begin{eqnarray}\label{13}
t=\frac{3s}{2}. 
\end{eqnarray}
For a higher invariant we get $t < \frac{3s}{2}$.

The relative quadratic helicity is well-defined. In a particular case the relative quadratic helicity express the Total 
invariant of magnetic braids, introduced in \cite{Y-H}.

\subsection*{An ergodic integral} 

Let us recall an approach by V.I.Arnol'd toward a description of the magnetic helicity integral as an asymptotic ergodic Hopf invariant of magnetic lines, see 
 \cite{A-Kh}. Assume that the magnetic field  $\B$ is represented by a finite collection of magnetic tubes $\Omega \subset \R^3$. Denote by 
$F(t): \Omega \to \Omega$ the ergodic magnetic flow along magnetic lines. By the Birkhoff Theorem the function  $(\A,\B)$ admits
the ergodic average, denote this function by $f: \Omega \to \R$. The function $f$ is invariant
with respect to the magnetic flow $F$. The function $F$ is constant restricted to each magnetic line - a trajectory of the flow $F$. The function $f$ is called the helicity density. 

The function $f$ is frozen-in with respect to volume-preserved transformations of $\R^3$. The function $f$ is integrable and, moreover, an arbitrary positive power $f^k = f \cdot \dots  \cdot f$, $k \in \Z, k \ge 1$ is integrable.
The magnetic helicity integral is the result of the integration of $f$ over $\Omega$. The quadratic magnetic helicity integral is the result of the integration $f^2$ over $\Omega$. 

The ergodic theorem allows us to transform a combinatorial invariant of closed magnetic lines into asymptotic invariants of magnetic fields with non-closed magnetic lines, details of the construction are in \cite{A}.
Invariants of magnetic fields are given by ergodic integrals.

\section*{A higher invariant of magnetic lines}

\subsection*{The $M$-Invariant for a triple of magnetic tubes} 

Consider a magnetic field
 $\B = \cup_i \B_i$ with a support into 3 magnetic tubes  $\Omega_i, \quad i=1,2,3$ correspondingly.
 Assume that inside the each magnetic tube a coordinate system 
 $U_i \cong D^2 \times S^1$ is fixed. Assume that this coordinate system corresponds with the standard volume form
 in $\R^3$ and the magnetic field $\B_i$ points strictly  along the $S^1$--coordinate of the system. 
This assumption simplifies calculations and gives no loss of a generality. 

The integral magnetic flow of 
$\B_i$ trough the cross-section of the magnetic tube  $U_i$ is denoted by  $\Phi_i$.
The integral linking number  $\int_{U_i} (\A_j,\B_i) dx = \Phi_i \Phi_j lk(i,j)$ of magnetic tubes 
$U_i$, $U_j$ is denoted by  $(i,j)$, $i,j = 1,2,3 \quad i \ne j$. 

A multivalued function with the period $(i,j)$, which is a restriction of the scalar branch of the vector-potential 
$\A_j$ on the magnetic tube $U_i$ 
denote by  $\varphi_{j,i}: U_i \to \R$. The function $\varphi_{j,i}$ is well-defined up to an additive constant. 
Consider a function 
$$ \phi_1 = (3,1)\varphi_{2,1} - (1,2)\varphi_{3,1}: U_1 \to \R, $$
which is well-defined by means of multivalued functions  
 $\varphi_{2,1}, \varphi_{3,1}$ up to an additive constant. To fix the constant, we assume that the following equation is satisfied:  $\int_{U_1} \phi_1 dx = 0$. Define the functions  $\phi_2$, $\phi_3$ by analogous formula. 

Define the vector
$$\F = (1, 3)(2, 3)\A_1 \times \A_2 + (2, 1)(3, 1)\A_2 \times \A_3 + (3, 2)(1, 2)\A_3 \times \A_1$$
$$ -\phi_1\B_1(2, 3) 
-\phi_2 \B_2(3, 1) -\phi_3\B_3(1, 2). $$
Obviously, the equation $\div(\F)=0$ is satisfied.

The vector-potential  $\G$, $\rot \G = \F$ and the integral  $\int_{\R^3} (\G,\F) dx$ are well-defined.
This integral is modified into the required invariant of volume-preserved diffeomorphisms. 
This modification includes the following extra 10 terms:

\begin{eqnarray}\label{e}
e_{1,2,3}=-2(1,2)(2,3)(3,1)(\int_{\R^3} \langle \A_1,\A_2,\A_3 \rangle dx)^2,
\end{eqnarray}

\begin{eqnarray}\label{14.1}
f_{1} = -2(\int_{U_1} \varphi_{2,1}^{var} (\grad
\varphi_{3,1}^{var},\B_1) dU_1)(\int_{\R^3} \langle \A_1,\A_2,\A_3 \rangle dx), 
\end{eqnarray}
\begin{eqnarray}\label{019}
d_{1,1} = -(2,3)^2 \int  \phi_1^2 (
\A_1,\B_1) dU_1,
\end{eqnarray}

\begin{eqnarray}\label{022}
d_{1;3} = (2,3)(1,2) \int  
\phi_1^2(\A_3,\B_1) dU_1. 
\end{eqnarray}

In the formula $(\ref{14.1})$ the terms $\varphi_{3,1}^{var}$, $\varphi_{2,1}^{var}$ are defined from 
$\varphi_{3,1}$, $\varphi_{2,1}$ correspondingly, see \cite{A2}. 
The extra 6 terms are defined by cyclic permutation of the indexes $\{1,2,3\}$ in the formulas $(\ref{14.1}),(\ref{019}),(\ref{022})$.

In \cite{A2} the following result is proved.

\begin{theorem}\label{1}
The integral expression
\begin{eqnarray}\label{M}
M(\B)=\int_{\R^3} (\G,\F) dx + e_{1,2,3}+ \sum_{i=1,2,3} f_{i}+d_{i,i}+d_{i;i+2} 
\end{eqnarray}
is invariant with respect to volume-preserved diffeomorphisms.
\end{theorem}


\subsection*{The invariant $M$ of 3 closed magnetic lines} 

Assume that  $\Phi_1=\Phi_2=\Phi_3=1$ and take a limit in the formula  $(\ref{M})$, when the thickness of
magnetic tubes tends to zero. The result satisfies the following definition.

\subsubsection*{Definition}
For an arbitrary 3-component link $\L \subset \R^3$
define a space (non-connected)
$Conf^{r}(\L) = (\L)^r$ 
as the Cartesian product of $r$ copies of $\L$.
The space $Conf^r(\L)$ is called the configuration space of the link $\L$.

Let
\begin{eqnarray}\label{conf}
F: Conf^r(\L) \to \R 
\end{eqnarray}
be an arbitrary integrable function on the configuration space.  
\[  \]

\begin{theorem}\label{conf}
Each term in the expression
 $(\ref{M})$ is defined by the integral of a corresponding function on the configuration space of the link $\L$.
\end{theorem}

Let us formulate an analogous definition for a function on the configuration space of magnetic lines.


\section*{Ergodic integrals and quasi-ergodic integrals.}

Let  
$\B$, $\div(\B)=0$ be a smooth magnetic field in  $\R^3$ with a support inside a finite collection of magnetic tubes 
$\Omega \subset \R^3$, the magnetic field $\B$ is tangent to the surface boundary of $\Omega$ and non-vanishes inside $\Omega$.

Define the configuration space
$K_{q,r}$ of magnetic lines, where  $p,r$ are given positive integers.  
For an arbitrary  $T>0$ (sufficiently great) consider collections of $r$ magnetic lines 
$L_1, \dots , L_r$ of  $\B$, the each line of the collection is parametrized by the standard segment  $[0,T]$ and is issued from points 
$\{l_1, \dots, l_r\}$ in  $\Omega$ correspondingly. Define collections of  $r(q+1)$ points, 
this collection consists of 
$r$ subcollections, $q+1$ points in each subcollection. The first subcollection  
$\{l_1; x_1, \dots, x_{q}\}$ consists of  $q+1$ points, including the initial point  $l_1$ on the magnetic line 
$L_1$. The second subcollection 
$\{l_2;x_{q+1}, \dots, x_{2q}\}$  consists of   $q+1$ points, each point is on the second magnetic line 
$L_2$, etc., the subcollection 
$\{l_r;x_{q(r-1)+1}, \dots, x_{qr}\}$ consists of  $q+1$ points, the each point belongs to the corresponding magnetic line точки, 
$L_r$. Obviously, the each point   $x_{qj+i}$ is well-defined by the corresponding parameter   $t_{qj+i}$, $1 \le j \le r$, $1 \le i \le q-1$, $0\le t_{qj+i} \le T$ of the magnetic flow, which transport the point 
$l_{j}$  to the point  $x_{qj+i}$ along the magnetic line $L_j$.  

 Let us say that the function  $F: K_{q,r} \to \R$ determines an ergodic integral, if
 the following conditions are satisfied:

$\bullet$ -1. For almost an arbitrary point  $\{l_1, \dots, l_r\} \in U^r$ 
the mean value $\bar{F}: K_{q,r} \to \R$ (in the sense of Ces\`aro) of the function $F$ with respect to 
position of points
$\{x_1, \dots, x_{q}, \dots, x_{q(r-1)+1}, \dots, x_{qr}\}$ is well-defined. 
By definition  $\bar{F}$ is induced by a function in the domain $U^r$ with respect to the projection  $\pi: K_{q,r} \to U^r$, $\pi(z)=\{l_1, \dots, l_r\}$, $z \in K_{q,r}$; denote this function by $\bar{F}: U^r \to \R$.

$\bullet$ -2. The function  $\bar{F}: U^r \to \R$ is locally integrable and is integrable.

The ergodic integral 
 $I(\B)$ is defined as the integral of the function  $\bar{F}$ over the domain $U^r$. 
\[  \]

Let us say that a function  
 $F: K_{q,r} \to \R$ determines a quasi-ergodic integral, if a linear mapping 
$X: K_{q,r} \to \R$ 
with respect to variables 
$\{x_1, \dots, x_{q}, \dots, x_{q(r-1)+1}, \dots, x_{qr}\}$ is well-defined, and, moreover,
for an arbitrary $p \in \R$ the restriction 
of $F$ to  $X^{-1}(p) \subset K_{q,r}$ satisfy Conditions 
 -1, -2; moreover, for an arbitrary  $p>0$ the integral  $f(p)=\int \bar{F} d(X^{-1}(p))$ determines 
an absolute bounded function 
$f(p): \R_+ \to \R, \quad p>0$. 
Additionally, if magnetic lines, issued from the points 
$\{l_1, \dots, l_r\}$ are closed, the function  $f: \R_+ \to \R$ is periodic.

The quasi-ergodic integral  $I(\B)$ is defined as a mean value of the function  $f(p)$ over 
$\R_+$. Generally speaking, this integral is multivalued and takes the value into a segment. 
In the case magnetic lines of $\B$ are closed, $f$ is periodic and a value $I(\B)$ is well-defined.

\begin{theorem}\label{MM}
The terms 
 $\int_{\R^3} (\G,\F) dx$, $e_{1,2,3}$,  $f_{i}$ in the formula  of  $M$, presented  in $(\ref{conf})$,
 are ergodic integrals. The terms 
$d_{i,i}$, $d_{i;i+2}$ in $(\ref{conf})$ are quasi-ergodic integrals. 
\end{theorem} 

In the paper \cite{A3} the following theorem is proved. 

\begin{theorem}\label{M}
Assume that magnetic lines of $\B$ inside $\Omega$ are closed. Then the invariant $M$ 
satisfy  Condition $C1$ for  $t=7$,  Condition $C2$ for  $s=12$, and Conditions $C3$, $C4$. 
\end{theorem}

\subsubsection*{Proof of Theorem  $\ref{MM}$}

A particular proof of Theorem is in  \cite{A4} (Theorem 3.1,(1) and Lemma 4.1.). 
I present a simplification of the proof for the main term
 $\int_{\R^3} (\G,\F) dx$ with simple estimations of the integral.
 
Assume that a magnetic field $\B$ is inside a finite collection of generic magnetic tubes,
$supp(\B) = \Omega \subset \R^3$.
 Recall the definition of the term $W$ of the integral $W$.  

Coordinates of a point in $K_{3,4;2}$ are given by collections 
$\{l_1,t_{1,1}, \dots t_{1,4},l_{2},t_{2,1}, \dots t_{2,4},l_3,t_{3,1}, \dots t_{3,4};y_1,y_2 \}$,
where $l_i \in U_i$, $t_{i,j} \in [0,T] \subset \R_{i,j}$, $j=1,2,3,4$, $y_1,y_2 \in \R^3$. 

Define the evolution mapping
 $F: K_{3,4;2} \to \Omega_1^{4} \times \Omega_2^4 \times \Omega_3^4$ by the formula 
 $$ F(l_1,t_{1,1}, \dots t_{1,4},l_{2},t_{2,1}, \dots t_{2,4},l_3,t_{3,1}, \dots t_{3,4} ) = $$
 $$(g^{t_{1,1}}(l_1), \dots g^{t_{1,4}}(l_1), g^{t_{2,1}}(l_2), \dots g^{t_{2,4}}(l_2), g^{t_{3,1}}(l_3), \dots g^{t_{3,4}}(l_3),$$
where $g^t$ is the magnetic flow of $\B$. From this formula the space 
 $K_{3,4;2}$ is the configuration space of $17$-points: $3(1+4)$ points $\{l_i, g^{t_{i,1}}(l_i), g^{t_{i,2}}(l_i), g^{t_{i,3}}(l_i), g^{t_{i,4}}(l_i)\}$
are on the magnetic lines, which are issued from $l_i$, $i=1,2,3$, and points $(y_1,y_2) \in (\R^3)^2$ are arbitrary.
 The standard volume form $dK_{3,4}$ on the space $K_{3,4;2}$ is well-defined.

The first step of the construction includes a definition of a function
 $W_{3,4;2}: K_{3,4;2} \to \R$, which is called the density function. 
The density function is not the lift of a function on  $\Omega^3$ by the projection  $\pi: K_{3,4;2} \to \Omega^3$.
The mean asymptotic value of the function  $W_{3,4;2}$ over the coordinates
 $t_{i,j}$, which is well-defined almost everywhere, depends of the parameters  $(l_1,l_2,l_3;y_1,y_2)$.
The last second step of the construction is a construction of a limiting tensor, this proves that the integral of $W_{3,4;2}$ over $\Omega^3 \times (\R^3)^2$ is well-defined.

Let us use the Gauss integral to calculate $W$ in the following formulas:
\begin{eqnarray}\label{171}
(2,3)(3,1)^2(1,2) \gamma_{t_{1,1},t_{2,1},t_{2,2},t_{3,1}}(\valpha_{1,2}(x_{1},x_{2,1};y_1),\valpha_{2,3}(x_{2,2},x_3;y_2)),
\end{eqnarray}
\begin{eqnarray}\label{172}
(2,3)^2(1,2)^2 \gamma_{t_{1,1},t_{1,2},t_{2,1},t_{2,2}}(\valpha_{1,2}(x_{1,1},x_{2,1};y_1),\valpha_{1,2}(x_{1,2},x_{2,2};y_2)).
\end{eqnarray}

In this formula by 
 $\gamma(\quad,\quad;\quad,\quad)$ is denoted the value of the kernel of the Gauss integral at a pair of corresponding vectors, the vectors of the pair depend of the parameters   $(x_1,x_{2,1},x_{2,2},x_3)$ and are attached
 to the points   $y_1,y_2$ correspondingly, the vectors $\valpha_{1,2}(x_{1},x_{2,1};y_1)$, 
 $\valpha_{2,3}(x_{2,2},x_3;y_2)$ in $(\ref{171})$ (for $(\ref{172})$ the formulas are similar) are given by $(\ref{R1})$, $(\ref{R2})$. 
The terms $(\ref{171})$, $(\ref{172})$ are well-defined in the asymptotic limit of all the positions of the points 
 $x_{i,j}$. For short we take $x_1=g^{t_{1,1}}(l_1)$, $x_3=g^{t_{3,1}}(l_3)$. The integration over
the variables $y_1,y_2$ is  taken after the asymptotic limit.

Let us investigate the term
 $(\ref{171})$ only, for the term $(\ref{172})$ the proof is analogous.  
The coordinates $\{t_{1,1}, \dots t_{1,4},t_{2,1}, \dots t_{2,4},t_{3,1}, \dots t_{3,4}\}$
are divided into the following 2 groups of coordinates, 
the coordinates of the first group  $\{t_{1,1},t_{2,1},t_{2,2},t_{3,1}\}$ are re-denoted by $\{\tau_{1}, \tau_{2,1}, \tau_{2,2}, \tau_3\}$ correspondingly. The coordinates of the second group  $\{t_{1,2},t_{1,3},t_{1,4},t_{2,3},t_{2,4},t_{3,2},t_{3,3},t_{3,4}\}$ are re-order as following:
$\{t_{1,2},t_{3,2},t_{1,4},t_{3,4},t_{1,4},t_{2,3},t_{2,3},t_{3,4}\}$ and
are re-denoted by
$\{\rho_{1,1},\rho_{3,1},\rho_{1,2},\rho_{3,2},\rho_{1,3},\rho_{2,3},\rho_{2,4},\rho_{3,4}\}$  correspondingly.

Let us define the factors in the formula $(\ref{171})$. 
Using the $4$ points of the first group $x_{1}=g^{\tau_{1}}(l_1), x_{2,1}=g^{\tau_{2,1}}(l_2), x_{2,2}=g^{\tau_{2,2}}(l_2), x_{3}=g^{\tau_{3}}(l_3)$, define the integral kernel  
\begin{eqnarray}\label{jadro}
\gamma_{\tau_{1},\tau_{2,1},\tau_{2,2},\tau_3}(\valpha_{1,2}(x_{1},x_{2,1}),\valpha_{2,3}(x_{2,2},x_3);y_1,y_2).
\end{eqnarray}  
Using the last $8$ points of the second group
$g^{\rho_{1,1}}=z_{1,1}, g^{\rho_{3,1}}=z_{3,1}, g^{\rho_{1,2}}=z_{1,2}, g^{\rho_{3,2}}=z_{3,2}, g^{\rho_{1,3}}=z_{1,3}, g^{\rho_{2,3}}=z_{2,3}, g^{\rho_{2,4}}= z_{2,4}, g^{\rho_{3,4}}=z_{3,4}$, define the integral kernel to calculate 
 $(2,3)(3,1)^2(1,2)$ by obvious way, see \cite{A-Kh} for the integral formula of the linking number.
The product of the expressions gives $(\ref{171})$. 

Let us prove that for almost arbitrary collection
$(l_1,l_2,l_3;y_1,y_2)$ there exists the asymptotic mean value of the expression
$(\ref{jadro})$ with respect to the variables
$\{\tau_{1}, \tau_{2,1}, \tau_{2,2}, \tau_{3}\}$. Denote this asymptotic mean value by
\begin{eqnarray}\label{meangamma}
\bar{\gamma}(l_1,l_2,l_3;y_1,y_2)
\end{eqnarray}


The absolute value of coordinates of the vector-potential 
 $\A(x_i;y)$ at an arbitrary point  $x_i \in L$ is integrable with respect to the parameter  $y \in \R^3$. 
This vector-potential determines the vector-functions  
\begin{eqnarray}\label{R1}
\valpha_{1,2}(x_1,x_{2,1};y_1) = \A(x_1;y_1) \times \A(x_{2,1};y_1), 
\end{eqnarray}
\begin{eqnarray}\label{R2}
\valpha_{2,3}(x_3,x_{2,2};y_2) = \A(x_3;y_2) \times \A(x_{2,2};y_2).
\end{eqnarray}
This vector-functions for  arbitrary fixed $y_1,y_2$ are integrable with respect to the parameters $\{x_{1},x_{2,1},x_{2,2},x_3\}$. 
 
 By the Birkhoff Theorem
 the vector-functions  $\valpha_{1,2}$, $\valpha_{2,3}$ in   $(\ref{jadro})$ 
admit the asymptotic limits with respect to the first group coordinates. 
The mean vector-functions are denoted by 
$\bar{\valpha}_{1,2}(l_1,l_2)(y_1)$, $\bar {\valpha}_{2,3}(l_2,l_3)(y_2)$, this vector-functions depend
formally of the points $(l_1,l_2,l_3)$, but, in fact, depend of the triple of magnetic lines $L_1,L_2,L_3$ only. 

The integral kernel 
 $(\ref{jadro})$ is calculated algebraically and  the term $(\ref{meangamma})$ is well-defined for almost arbitrary $(l_1,l_2,l_3;y_1,y_2)$. 
Analogously, the integral kernel  $W_{3,4;2}$, corresponded to $(\ref{171})$,  admits the mean value 
over all the variables $\{\tau_{1}, \tau_{2,1}, \tau_{2,2}, \tau_3;\rho_{1,1},\rho_{3,1},\rho_{1,2},\rho_{3,2},\rho_{1,3},\rho_{2,3},\rho_{2,4},\rho_{3,4}\}$.
Denote this mean value by
\begin{eqnarray}\label{meanW}
\bar{W}(l_1,l_2,l_3;y_1,y_2).
\end{eqnarray}
 
 The product of pairwise asymptotic linking numbers is well-defined for almost arbitrary
collections of pairs of magnetic lines $(l_1,l_2), (l_2,l_3), (l_3,l_1)$  (see \cite{A}).  The vector 
 $(\ref{meanW})$ is well-defined and
  the first step of  the construction is described.
 
  Pass to the second step of the construction and prove that the integrals $(\ref{meangamma})$,
$(\ref{meanW})$   over $\Omega^3 \times (\R^3)^2$ are well-defined.
   Estimate the total term $(\ref{171})$
by a limiting tensor, which is absolutely integrable over the configuration space 
 $\Omega^3 \times (\R^3)^2$. Denote by  $a(x_1,x_{2,1},x_{2,2},x_3)$ 
the absolute value of the term 
 $(\ref{jadro})$ (the value $+\infty$ is admitted) after the integration over the variables  $y_1, y_2$. 
Assume firstly that the points $x_1,x_{2,1},x_{2,2},x_3$ belong to the triple of the segments of magnetic lines,
which are pairwise close to each other. Denote by $\delta$ a small parameter, which   
is the distance of the segment on $L_2$ to the segments on  
 $\{L_1,L_3\}$ (for short we assume that the segment on $L_1$ is closer that the segment on $L_3$ to the segment on $L_2$). 

\begin{lemma}$\label{77}$
Let $y_1=y_2=x_{1,1}=x_{1,2}=x_2=x_3$, and $\omega > 0$ be a given positive (arbitrary small) number,
constants $\delta_0 > \delta$ be arbitrary.
Take an arbitrary non-degenerate $\delta$-variation of the magnetic line  $L_1$ and  
an arbitrary variation of the magnetic line $L_3$, which is estimated from above by $\delta$ and from below
by  $\delta_0$. Then the absolute integral value of the term 
\begin{eqnarray}\label{avarepsilon}
a(x_1,x_{2,1},x_{2,2},x_3)
\end{eqnarray}
over arbitrary 
$\varepsilon$--variations of points $y_1,y_2,x_1,x_{2,2},x_3$ (the point $x_{1,1}$ is fixed)
along the corresponding segments of magnetic lines  is estimated by 
$C \delta^{-1-\omega}$, where the positive constant $C$ depends only on $\varepsilon$.
The constant $\varepsilon$ depends on the norm of the 2-jets of $\B$ in $\Omega$ and depends no of  
$\delta$.  
\end{lemma}

\subsubsection*{Remark}
By the results of \cite{A4}, one may replace $C \delta^{-1-\omega}$ by  $C\log(\delta^{-1})$ in the lemma. 

\subsubsection*{Proof of Lemma  $\ref{77}$}
 
To simplify the notation put 
$\varepsilon=1$. The singularity in the configuration space is of the order 
 $r^{-10}$, where  $r$ is the distance in $\R^3$ which corresponds to the parameter of deformation. 
This formal order includes the order $-2$ of the each magnetic dipole (4 dipoles), the order of the kernel in the Gauss integral, given by    
  $\dist(y_1,y_2)^{-2}$. The integration of the term  $(\ref{avarepsilon})$ is over the 6-dimensional domain
of the variables   
 $y_1,y_2$ and of a 3-dimensional domain, of the variables   $x_{1,1},x_{1,2},x_2,x_3$. 
 As the result, we get that the singularity  of $(\ref{avarepsilon})$ is of the formal order
 $-1$. 

After the deformation, described in the lemma, the term $(\ref{avarepsilon})$ is well-defined and integrable. 
To calculate this generic term, we integrate singular functions of the order
 $r^{-6}$ (the coordinate $r$ is the distance between the parameters  $x_{1,1}, x_{1,2}$ on the line  $L_1$)
over 7-dimensional space.  The integral is well-defined.
A formal estimation  (over the parameter $\delta$) of the deformation of the singularity proves Lemma.  
$\ref{77}$.

\[  \]

Let us estimate  $W$ by absolute value using the lemma. Consider the cube with the edge of the length  $T$ in the configuration space, which is given by the parameter of the magnetic flow. The configuration space is a union of a finite number of small cubes. Let us define a limiting tensor of  
$W_{3,4;2}$ in each cube. Recall that the  limiting tensor is absolutely integrable over the configuration space.
and estimates 
the absolute value of $W_{3,4;2}$.

We start with cubes, called diagonal cubes, which are closed to top singularities, which are described in Lemma $\ref{77}$, up to parallel translations of all 4-points along the magnetic flow.  The last cubes in the configuration space, called peripheral cubes,  are defined analogously. 

In each diagonal cube we get the estimation from Lemma $\ref{77}$. In an arbitrary peripheral cube 
estimations is more simple, and formally are given by the same formulas, $\delta_0$ is not a small parameter. 
As the result we get that the expression 
$(\ref{jadro})$ is estimated by a function of the order $\delta^{-1-\omega}$, where $\delta$ is the minimal
pairwise distance between segments of magnetic lines (if there is a pair of close segments of magnetic line) and by a function of the order $1$,
if all the segments are pairwise non-closed.

By the Holder inequality we get:
$$ \int fg dx \le (\int \left|f\right|^q dx)^{\frac{1}{q}} (\int \left|g\right|^p dx)^{\frac{1}{p}}, \quad 1<p<2, \quad \frac{1}{p} + \frac{1}{q}=1. $$
In this inequality $f$ is the limiting tensor for 
 $(\ref{jadro})$,  $g$ is the limiting tensor with logarithmic singularities for the term
 $(2,3)(3,1)^2(1,2)$, which is much simple. We use the denominator  
 $p=1+\omega$ and a large denominator  $q$. The function $(\ref{meanW})$ is integrable and the main term $W$ is given by an ergodic integral.

\section*{Examples of Magnetic Knots in the standard sphere $S^3$ and in several homogeneous manifolds}

In this section we consider examples of magnetic knots with closed magnetic lines (or with magnetic lines on family of surfaces) inside  compact (homogeneous) manifolds,
for which $M$-invariant is  non-vanished. The Examples $I$ and $II$ are generalizations with non-simply connected manifolds. For this examples Theorems \ref{1},\ref{conf},\ref{M} are conjectured. 

\subsection*{A one-parametric family of magnetic knots in  $S^3$ }

Consider the standard singular fibration $S^3 \to S^2$ with 2 singular linked circles $S^1_1 \subset S^3$, $S^1_2 \subset S^3$, and
with Hopf family of regular tori $T_t$, $t \in [1,2]$ between this two circles, $T_1, T_2$ are shrined into
$S^1_1$ and $S^1_2$ correspondingly.  Consider the Cartesian coordinate system $(x,y,z)$ on $S^3 \setminus \{\infty\}$.
The circle $S^1_1$ is the unite central circle on the plane $(x,y)$. The circle $S^1_2$ is the standard vertical $z$-axis, $\infty \in S^1_2$, through the origin.

Define a real parameter $r$,  $1\le r \le 2$.  Define a $r$-parameter family  of magnetic knots $\Upsilon_r$ in $S^3$.
Magnetic lines of $\Upsilon_r$ for each $t \in ]1,2[$ are on $T_t$ and wind
$1$ time along the $S_1^1$--parallel of $T_t$ and $r$ times along the $S^1_2$-meridian of $T_t$. For rational $r$,
the magnetic knot $\Upsilon_r$ consists of closed lines. The magnetic knots $\Upsilon_{r_1}$, $\Upsilon_{r_2}$, in the case
$r_1 \ne r_2$, are not equivalent with respect to volume-preserved diffeomorphisms of $S^3$. In the case  $r=1$ 
we get the standard Hopf fibration with fibers along the standard Hopf mapping $h: S^3 \to S^2$.

The combinatorial formula of the invariant $M(\L)$, in the case $\L$ is a 3-component link, is well-defined up to
the sum with a polynomial $P((1,2),(2,3),(3,1))$, which depends on pairwise linking numbers  of $\L$
(see \cite{A3}). In the case $(1,2)=(2,3)=(3,1)$, to keep asymptotic properties of $M$, we  assume that 
$\deg(P(k)) \le 11$. We define $\dot{M} = M + P$, the invariant $\dot{M}$ is ergodic. Moreover,
without lost of a generality we assume that $\dot{M}$ is trivial on a prescribed collection of the following 2 simplest 
links $\L_{1,1,1}$, $\L_{2,2,2}$.

The link $\L_{1,1,1}$ consists of 3 magnetic lines with pairwise linking number $1$,
each line is a fiber of the Hopf fibration $h: S^3 \to S^2$. By the construction,
$\L(1)=\L_{1,1,1}$, where $\L(1)$ is the link, which is defined by an arbitrary ordered magnetic lines of the
magnetic knot $\Upsilon_1$. 

The link $\L_{2,2,2}$ is defined as following.
Take the symmetric triangle with the unite edges on the plane. Take 3 circles $L'_1,L'_2,L'_3$ of the radius $\frac{1}{2}$
around its vertexes, which are tangent to each other in the centers of edges. Then take a small 3D deformation
of $(L'_1,L'_2,L'_3) \to (L_1,L_2,L_3)$ in small neighborhoods of tangent points of the pairs $(L'_1,L'_2)$, $(L'_2,L'_3)$, $(L'_3,L'_1)$; as the result we assume that the pairwise linking numbers
of $(L_1,L_2)$, $(L_2,L_3)$, $(L_3,L_1)$ are equal to $+2$. 

Denote by $\L(2)$ a 3-component link, which is defined by an arbitrary ordered triple of generic magnetic lines of the
magnetic knot $\Upsilon_2$
It is not difficult to prove that  pairwise linking numbers of $\L(2)$
and   $\L_{2,2,2}$ coincide.  

By the construction $\L_{2,2,2}$ is distinguished from  $\L(2)$ by the commutator of 3-components (or, equivalently, by the $\Delta$-moves of 3 components). By the following lemma  and the combinatorial formula of $M$ from [\cite{A3}(17)], the value $M(\L(2))$ is distinguished from   $M(\L_{2,2,2})$
by a non-zero integer.

\begin{lemma}\label{semi-boundary}

Let  $\L = (L_1 \cup L_2 \cup L_3)$ be an arbitrary 3-component link for which the pairwise linking coefficients
$(1,2), (2,3), (3,1)$ are even. Let $\L' = (L'_1 \cup L'_2 \cup L'_3)$ be the 3-component link, which is the result of a $\Delta$--move of $\L$ with 3 different components.

The parity of the coefficients $C_2(\L)$, $C_2(\L')$ of the Conway polynomial are distinguished,
 and the invariants 
$Arf(\L)$, $Arf(\L')$ are distinguished.  
\end{lemma}

\begin{remark}
The invariant 
$Arf(\L)$ is well-defined in a kess restrictive case, when  all the pairwise linking numbers of $\L$ are odd. 
\end{remark}

\subsection*{Proof of Lemma  $\ref{semi-boundary}$}
For a link $\L = (L_1 \cup L_2 \cup L_3)$ which satisfies the lemma, the equation 
 $\mu^2_{123}(\L) \equiv C_2(\L) \pmod{GCD(1,2),(2,3),(3,1)}$ is proved in  \cite{M}, Theorem 3.5.
The equation $Arf(\L) \equiv \mu_{123}(\L) \pmod{2}$ is proved using the Gauss diagrams as in \cite{M-P}.
 $Arf$-invariant satisfy the lemma. 
Lemma $\ref{semi-boundary}$ is proved.
\[  \]

The invariant $M$ is ergodic, therefore $M(\L(r))$ in continuously changed from $M(\L(1))=0$ to
$M(\L(2)) \ne 0$, $1 \le r \le 2$.

\subsection*{Example  ${\rm{I}}$ of a magnetic knot in the rational homological sphere   $S^3/\Q$}

In the group of unit quaternions $S\H$ consider the subgroup of integer quaternions  $\Q \subset S\H$ 
$$\{ \i,\j,\k \quad \vert \quad
\i\j = \k = -\j\i,
\j\k=\i = -\k\j, \k\i = \j = -\i\k, \i^2 =\j^2 = \k^2 =-1 \}.$$ 

Consider the standard (right) action 
$\Q \times S^3 \to S^3$,  which is well-defined because of the diffeomorphism $S\H \cong S^3$. 
Consider the 2-sheeted covering
 $S\H \to SO(3)$, the image of the subgroup  $\Q \subset S\H$ is the Klein subgroup  $\K \subset SO(3)$, $\K \cong \Z/2 \times \Z/2$. The Klein group acts on 
$S^2$, this action is induced by the standard projection $SO(3) \to S^2$, the action has 6 fixed points, which are 
the intersection points of the standard unite sphere $S^2 \subset \R^3$ with the coordinate axis. The elements of $\K$
acts on $S^2$ by rotations trough the angle $\pi$ with respect to the corresponding coordinate axis.

The following commutative diagram of groups 
\begin{eqnarray}\label{Q}
\begin{array}{ccc}
\Q \times S^3 & \to & S^3/\Q \\
\downarrow & & \downarrow \\
\K \times S^2 & \to & S^2/\K, 
\end{array}
\end{eqnarray}
is well-defined. In this diagram horizontal maps are projections onto the orbits of the action,
the left vertical mapping is the Cartesian product of the projection 
  $\Q \to \K$ and the composition 
$S^3 \cong S\H \to SO(3) \to S^2$, which coincides with the standard Hopf fibration,
the right vertical mapping is induced from the left vertical mapping by the projection onto the orbits.

The magnetic knot in $S^3/\Q$ with closed magnetic lines is well-defined by fibers of the right vertical mapping in the diagram. 
A generic magnetic line $L \subset S^3/\Q$ of this magnetic knot represents an oriented cycle  $[L] \in H_1(S^3/\Q;\Z)$, which is not an oriented boundary, but is a non-oriented boundary. This means that the magnetic line is a boundary of a non-oriented  Seifert surface.  The manifold   $S^3/\Q$ admits a natural trivialization of the tangent bundle,
$T(S^3/\Q) \cong 3\varepsilon$. Take a magnetic line $L_1$ in the integer homology class of $[L]$. From this data the Arf-Brown invariant  $\Theta(L_1) \in \Z \pmod{8}$ is well-defined.

\subsubsection*{A generalization of the Example}

The diagram $(\ref{Q})$ is included into the following diagram: 
\begin{eqnarray}\label{I} 
\begin{array}{ccc}
\Sigma \times S^3 & \to & S^3/\Sigma \\
\downarrow & & \downarrow \\
\I \times S^2 & \to & S^2/\I. 
\end{array}
\end{eqnarray}
In this diagram $\Q \subset \Sigma$ is the Poincar\'e extension of the index $15$ of the integer quaternions
to the fundamental group of the integer homology sphere, 
 $\K \subset \I$ is the extension of the Klein group to the icosahedron group, the lower horisonatal mapping of the diagram is a free action,  the bottom mapping is the semi-free action.  By the Klein uniformization
 \cite{K}, the icosahedron group  $\I$ is covered by the modular group  $PSL(2,\Z)$,
 which acts conform on the half-plane.

Below in the diagram  $(\ref{seq2})$ a quadratic extension 
$\Q \subset \aleph$ is well-defined. The quadratic extension $\Q \subset \aleph$ is mapped into a quadratic 
extension
$\K \subset \D$ by the projection onto the factorgroup, where  $\D$ is the dihedral group of the order $8$.

The inclusion 
 $\K \subset \I$ admits no  extension of the quadratic extension 
$\K \subset \D$  of the subgroup to a quadratic extension of the group $\I$. 
The minimal infinite-order extension 
 $\I \subset \Upsilon$ is well-defined, where $\Upsilon$ is a Kleinian group, which acts conform on $\hat{\C}$,
 and this action extends the action of the Fuchsian group. The group $\Upsilon$ is covered by a group, which acts conform on the half-plane.

\subsection*{Example ${\rm{II}}$ of magnetic knot in the rational homological sphere  $S^3/\Q$}
The standard Hopf fibration 
 $h: S^3 \to S^2$, is given by the formula  
$ \{ (z_1, z_2) \}, \vert z_1 \vert^2 + \vert z_2 \vert^2 = 1 $,
$$h:(z_1,z_2) \mapsto \frac{z_1}{z_2}.$$
The conjugated Hopf fibration 
$\bar{h}: S^3 \to S^2$ is given by the formula 
$$\bar{h}:(z_1,z_2) \mapsto \frac{\bar{z}_1}{z_2}.$$
The following diagram  
$$ 
\begin{array}{ccc}
\Q \times S^3 & \to & S^3/\Q \\
\downarrow & & \downarrow \\
\Z/2 \times S^2 & \to & \RP^2, 
\end{array}
$$
is well-defined, where $\Q \to \Z/2$ is the epimorphism with
the generator $\i$ is the kernel, 
$\Z/2 \times S^2  \to  \RP^2$ is the projection of the antipodal involution, see  \cite{S}.

Define the magnetic knot on $S^3/\Q$ by the fibers of 
 $\bar{h}$. An arbitrary magnetic line $L \subset S^3/\Q$ of the magnetic knot is not an non-oriented boundary. 
The group
 $\Q$, which is the fundamental group of the rational homology sphere  $S^3/\Q$ admits a quadratic extension 
 $\Q \subset \aleph$, which is defined below by $(\ref{seq2})$. By this extension the image of the generator $\i \in \Q$ in $\aleph$ belongs to the commutant $[\aleph,\aleph] \subset \aleph$. 
A Seifert surface for $L$ is well-defined as a surface with a prescribed normal bundle structure  (see below  
the definition of this structure in Theorem \ref{th2}) with a control to the Eilenberg-MacLane space  $K(\aleph,1)$.
For Seifert surfaces with prescribed normal bungle structures the hyperquaternionic Arf-invariant is well-defined as an integer 
 $\pmod{16}$.

Examples 
 $I$, $II$ of magnetic knots on  $S^3/\Q$ assume that asymptotic ergodic $M$-invariant is generalized for 
magnetic knots in rational homology spheres. The parity of $C_2$-coefficient of the Conway polynomial for classical links in $\R^3$ corresponds to   
the Arf-invariant.  In the next section  we determine a group $W$, which is called the Witt group of hyperquaternionic forms.  
The reason to introduce the hyperquaternionic Arf-invariant is clarify by the following diagram: 
$$ 
\begin{array}{ccc}
 C_2 \ of \ the \ Conway \ polynomial & \longrightarrow & Arf \ invariant \\
 &  &  of \  classical \ links\\
\downarrow &   &  \downarrow \\
? & \longrightarrow &  hyperquaternionic \ Arf-invariant \\
&  &  of \ links \ in \ S^3/\Q.\\
\end{array}
$$
In the diagram by 
 $?$ is denoted a hypothetic integer-valued finite-type invariant of links in rational homological spheres,
 which determines asymptotic ergodic invariants.

\section*{Hyperquaternionic \ Arf-invariant}

\subsection*{Arf-invariants of immersed surfaces} 

Consider an immersion  $\varphi: M^2 \looparrowright \R^3$ of a closed, generally speaking, non-oriented surface
into $\R^3$.  The immersion  $\varphi$ up to regular cobordism  represents an element of the group denoted by $Imm^{sf}(2,1)$, we use notations as in \cite{A-E}. The Arf-Brown invariant is an isomorphism  $$\Theta: Imm^{sf}(2,1) \to \Z/8.$$ 
Denote 
$Imm^{sf}(2,1)$ by $V$ for short (an algebraic definition of $\Theta: V \cong \Z/8$, using $\Z/4$-quadratic forms, is in  \cite{G-M}).
If $M^2$ is an orientable surface, the element $\Theta([\varphi])$ belongs to the subgroup $\Z/2 \subset \Z/8$.
In this case the element $\frac{\Theta([\varphi])}{4} \pmod{2}$ is called the Arf-invariant of $[\varphi]$.

Let  $K^3$ be a closed oriented 3-dimensional manifold. Assume that a trivialization of the tangent bundle
$\Psi: T(K^3) \cong 3\varepsilon$ is fixed. Assume that an immersion  $\varphi: M^2 \looparrowright \R^3$ of a closed surface is given. The immersion 
 $\varphi$ represents an element $[\varphi]$ in the group  $V$ and the Arf-Brown invariant  $\Theta([\varphi])$ is well-defined. 

In the case, when $M^2$ is a surface with a boundary, assume that each component of the immersed curve $\varphi(\partial M^2)$ has the trivial stable Hopf invariant (= an even self-linking number). In this case 
the Arf-Brown invariant $\Theta([\varphi])$ is well-defined.

\subsection*{Group $\aleph$ of the order $16$}

Consider the cyclic group $C_8$ of the order $8$,   
 $C_8=\{\exp{(\frac{k\pi \i}{4})} \quad \vert k \in \Z/8\}$.
 Denote by 
$\theta: C_8 \to C_8$, $\theta: S \mapsto S^3$, $S \in C_8$ the cubing automorphism. 
Let us define a group
$\aleph$ of the order $16$, by attaching an element  $T$ of the order $2$ by the equation  $TST=S^3$, 
see for details
 \cite{C-M}, Ch.1 1.8. The following short exact sequence:
\begin{eqnarray}\label{seq1}
 0 \to C_8 \to \aleph \to \Z/2 \to 0
\end{eqnarray} 
is well-defined. In this sequence the left mapping is  the inclusion on the subgroup, the right mapping is the projection onto. 
Denote by  $T\j \in C_8$ a generator of the subgroup 
$C_8 \subset \aleph$; denote by $\j \in \aleph$ the element   
$T(T\j)$; denote by $\k \in \aleph$ the element $T\j T$; denote by $-1 \in C_8 \subset \aleph$ the  element 
 $(T\j)^4$, denote by $-\i$ the element 
 $(T\j)^2=\k\j$.

Define the following short exact sequence
\begin{eqnarray}\label{seq2}
 0 \to \Q \to \aleph \to \Z/2 \to 0, 
 \end{eqnarray} 
 where $\Q$ is the integer quaternions subgroup.
The group  $\Q$ is of the order $8$, this group admits the following 
standard corepresentation:
$$\{ \i,\j,\k \quad \vert \quad
\i\j = \k = -\j\i,
\j\k=\i = -\k\j, \k\i = \j = -\i\k, \i^2 =\j^2 = \k^2 =-1 \},$$ 
which corresponds to the notations of the generators.

\subsection*{Representation $\Phi: \aleph \to \SO(4)$}
Define a  $\SO(4)$--representation 
$\Phi: \aleph \to \SO(4)$ by the following matrices: 
\begin{eqnarray}\label{repr}
\Phi(T)= 
 \begin{array}{cccc}
1 & 0  &  0  &  0 \\
0 & -1 &  0  &   0 \\
0 &  0  &  0  &  1 \\
0  &  0  &  1  &  0 \\
\end{array}
\qquad \qquad \Phi(T\j)= 
 \begin{array}{cccc}
0 & 0  &  -1  &  0 \\
0 & 0 &  0  &   -1 \\
0 &  -1  &  0  &  0 \\
1  &  0  &  0  &  0 \\
\end{array}
\end{eqnarray}

The elements
 $\i$, $\j$, $\k$ are given by the following matrices: 
\begin{eqnarray}\label{repr1}
\Phi(\i)= 
 \begin{array}{cccc}
0 & -1  &  0  &  0 \\
1 & 0 &  0  &   0 \\
0 &  0  &  0  &  -1 \\
0  &  0  &  1  &  0 \\
\end{array}
\qquad \Phi(\j)= 
 \begin{array}{cccc}
0 & 0  &  -1  &  0 \\
0 & 0 &  0  &   1 \\
1 &  0  &  0  &  0 \\
0  &  -1  &  0  &  0 \\
\end{array}
\qquad \Phi(\k)= 
 \begin{array}{cccc}
0 & 0  &  0  &  -1 \\
0 & 0 &  -1  &   0 \\
0 &  1  &  0  &  0 \\
1  &  0  &  0  &  0 \\
\end{array}
\end{eqnarray}

The representation 
 $\phi = \Phi \vert_{\Q}: \Q \to \SO(4)$  is equivalent to
the standard representation  
 $\Q \to \SH \subset \SO(4)$.

\subsection*{The octahedral extension  $\aleph \subset \Upsilon$ of the index $3$}
Let us unify short exact sequences 
 $(\ref{seq1})$, $(\ref{seq2})$ into the following diagram: 
\begin{eqnarray}\label{seq2D}
\begin{array}{ccccccccc}
  &     & 0                 &         & 0 &     &      &     & \\
  &     & \uparrow                 &         & \uparrow &     &      &     & \\
0 & \to & \K & \subset & \D & \to & \Z/2 & \to & 0 \\
  &     & \uparrow                &         & \uparrow &     &  \|    &     & \\
0 & \to & \Q  & \subset & \aleph & \to & \Z/2 & \to & 0 \\
  &     & \uparrow                  &         & \uparrow &     &      &     & \\
    &     & \Z/2                 &     \cong    & \Z/2 &     &      &     & \\
  &     & \uparrow                 &         & \uparrow &     &      &     & \\
    &     & 0                 &         & 0 &     &      &     & \\    
\end{array}
\end{eqnarray}

In this diagram by $\D$ is denoted the dihedral group of the order $8$, the projection 
 $\aleph \to \D$ extends the reduction $C_8 \to C_4$ of the cyclic subgroup modulo 4, $\K \cong \Z/2 \times \Z/2 \subset \D$ is the Kleinian group, $\Q \to \K$ is the natural epimorphism, which  is the projection onto the  central quotient   $\{\pm 1 \} \subset \Q$. The group $\K$ is equipped with the representation
 $\K \to \SO(3)$, the image of the corresponding element $[\i]$, $[\j]$, $[\k]$ 
is the rotation trough the angle $\pi$ with respect to the axis, which is perpendicular to the coordinate plane  
$P_{\i}, P_{\j}, P_{\k}$ in $\R^3$ correspondingly. 

The group 
$\D$ is equipped with the representation  $\tilde \lambda: \D \to \O(3)$, the element $[T] \in \D$, which is define as the image of the element $T \in \aleph$ by the projection
$\aleph \to \D$,  is represented by  symmetry with respect to the plane, which is perpendicular to
 $P_{\i}$, along the bisector of the coordinate planes   $P_{\j}$ and $P_{\k}$.  
The representation  $\tilde \lambda \vert_{\K} = \lambda$ is defined such that the representation  $\phi: \Q \to \S^3 \subset \SO(4)$ covers the representation $\lambda$ by the projection 
$S^3 \to \SO(3)$. 

The representation 
$\tilde \lambda: \aleph \to \O(4)$ is the quadratic extension of the representation 
$\lambda$ by the standard quadratic extension  $\SO(3) \subset \O(3)$.  

Define the following diagram: 
\begin{eqnarray}\label{seq3D}
\begin{array}{ccccccccc}
  &     & 0                 &         & 0 &     &      &     & \\
  &     & \downarrow                 &         & \downarrow &     &      &     & \\
0 & \to & \Z/3 \tilde{\times} \K & \subset & \Z/3 \tilde{\times} \D & \to & \Z/2 & \to & 0 \\
  &     & \downarrow                &         & \downarrow &     &  \vert \vert    &     & \\
0 & \to & \I  & \subset & \Upsilon & \to & \Z/2 & \to & 0 \\
\end{array}
\end{eqnarray}

The group  
$\Z/3 \tilde{\times} \K$ is a semi-direct product of the subgroups  $\D$, $\Z/3$ in the icosahedron group. 
The subgroup  $\Z/3$ permutes the images of quaternion units  $[\i]$, $[\j]$, $[\k]$ in $\K$. 
The inclusion 
 $\Z/3 \tilde{\times} \K \subset \I$ into the icosahedron group is of the index $5$. 
 The group 
$\Upsilon$ is the fundamental group of the  homology Poincar\'e sphere. The inclusion 
 $\Z/3 \tilde{\times} \D \subset \Upsilon$ is the quadratic extension of the inclusion 
$\Z/3 \tilde{\times} \K \subset \I$.

\begin{lemma}
--1. Diagram $(\ref{seq3D})$ is well-defined and contains the diagram  $(\ref{seq2D})$ as a subdiagram.

--2. The groups  $\Upsilon$, $\B$ are equipped with representations  $M:\Upsilon \to  \SO(4)$,
$\mu: \Z/3 \tilde{\times} \D \to \O(3)$, the representations $M, \mu$ extend the representations  $\Phi$, $\tilde \lambda$ correspondingly. 
\end{lemma}

\subsection*{The Witt group  $W$  of  hyperquaternionic forms}

Define the regular cobordism group of closed surfaces, the elements of $W$ will be called  hiperquaternionic forms.
Denote this group by $W$, from algebraic point of view, $W$ is a Witt group of special quadratic forms.  

Define an epimorphism  
 $\alpha: \aleph \to \Z/2 =\{\pm 1\}$ by the following formula:  $T\j, T \in \aleph$,  $\alpha(T\j)=-1, \alpha(T)=+1$. The kernel $Ker(\alpha)$ coincides with the dihedral subgroup  $\D \subset \aleph$.

Define an epimorphism 
 $\beta: \aleph \to \Z/2 =\{\pm 1\}$ by the following formula: $\beta(T\j)=-1, \beta(T)=-1$. The kernel $Ker(\beta)$  coincides with the quaternion subgroup $\Q \subset \aleph$.

Over the space
 $B\aleph=K(\aleph,1)$ the canonical vector $\SO(4)$--bundle is well-defined, 
the structure group of the canonical bundle is defined by the representation  
 $\Phi: \aleph \to \SO(4)$, denote this universal bundle by  $A$. Denote by $\gamma$ the line canonical bundle over 
 $B\Z/2 \cong \P^{\infty} \cong K(\Z/2,1)$. Denote by $\alpha: K(\aleph,1) \to K(\Z/2,1)$ the mapping of the classifing spaces, which is associated with  the homomorphism
$\alpha$, denote by $\beta: K(\aleph,1) \to K(\Z/2,1)$ the mapping, which is associated with the homomorphism $\beta$.

A triple $(M^2, \eta_M, \Xi_M)$ is called a hyperquaternionic form, where

$\bullet$ $M^2$ is a closed, generally speaking, non-orientable surface; 

$\bullet$ $\eta_M = : M^2 \to K(\aleph,1)$ is a characteristic class, the composition  $\alpha \circ \eta_M$ 
is denoted by
 $\eta_{\alpha;M}: M^2 \to K(\Z/2,1)$, the composition
 $\beta \circ \eta_M$ is denoted by  $\eta_{\beta;M}: M^2 \to K(\Z/2,1)$;

$\bullet$ $\Xi_M$ is the isomorphism  $T(M) \oplus \eta_{\alpha;M}(\gamma) \oplus \eta_M^{\ast}(A) \oplus 3\eta_{\beta;M}(\gamma) \cong 10\varepsilon$, where by 
$\varepsilon$ is the trivial line bundle.  

In particular, by definition of 
 $\Xi_M$, the characteristic class 
$\eta_{\alpha;M} + \eta_{\beta;M}: M^2 \to K(\Z/2,1)$ corresponds to the orientation homomorphism  $H_1(M;\Z/2) \to \Z/2$, (denote $\eta_{\alpha;M} + \eta_{\beta;M} = \kappa_M: M^2 \to K(\Z/2,1)$,  this characteristic class coincides with 
the characteristic Stiefel-Whitney class $w_1(M)$). 

On a set of all hyperquaternionic form an additive operation by a disjoint union is well defined. 
The standard regular cobordism relation determines an equivalence  relation of quadratic hyperquaternionic forms.
The cobordism group up to this equivalence relation is denoted by  
 $W$, this is the required Witt group.

\begin{definition}\label{hyper}
A hyperquaternionic form 
$(M^2, \eta_M, \Xi_M)$, for which the characteristic mapping $\eta$ takes values in the subspace 
 $K(\Q,1) \subset K(\aleph,1)$, is called a quaternionic form.
\end{definition}

\begin{theorem}\label{th2}
The group $W$ contains a cyclic subgroup  $P \subset W$ of the order $16$, $P \cong \Z/16$.


\end{theorem}

\begin{definition}
Define a subgroup $W_{\Q} \subset W$ in the Witt group as the group, which is generated by quaternionic forms. Define the group  $W^{\odot}_{\Q}$, which is called the Witt group of quaternion forms. The group  $W^{\odot}_{\Q}$
is generated by quaternion forms, the regular cobordism relation for this group 
assumes the following additional property: 

$\bullet$ the structure mapping on a cobordism manifold admits a prescribed reduction to a mapping with the image in the quaternion classifying subspace  $K(\Q,1) \subset K(\aleph,1)$.
\end{definition}

By the construction, the canonical projection 
 $p: W^{\odot}_{\Q} \to W_{\Q}$ is well-defined.

\subsection*{The Arf-Brown homomorphism the group  $W_{\Q}$ onto the Witt group of $\Z/4$--quadratic forms} 

Denote by 
 $V$ the Witt group of $\Z/4$--quadratic forms with Arf-Brown invariants. This group
is related with the Rokhlin's Signature Theorem, see 
\cite{G-M}. The group $V$ is the cyclic group of the order $8$. Define the forgetful homomorphism 
$$\rho^{\odot}: W_{\Q}^{\odot} \to V$$
from the Witt group of quaternionic forms into the Witt group of quadratic
$\Z/4$--forms as following.

Let 
$(M^2, \eta, \Xi)$ be a quaternionic form represented an element in  $W^{\odot}_{\Q}$. 
Consider the standard 3-skeleton 
$S^3/\Q \subset K(\Q,1)$, which is represented by the standard quaternion lens space. 
The pull-back of the bundle 
$A$ over $K(\aleph,1)$ with respect to the inclusion 
$S^3/\I \subset K(\I,1) \to K(\aleph,1)$ is denoted by  $A_{S^3/\Q}$.
The canonical isomorphism 
 $A_{S^3/\Q} \cong 4\varepsilon$ of the vector bundles over  $S^3/\Q$ is well-defined. 
The pull-back isomorphism 
 $\eta^{\ast}(A_{S^3/\Q}) \cong 4\varepsilon$, determines the isomorphism  $\Xi_M: \nu_M \to 7\varepsilon \oplus \kappa$, where $\nu_M$  is the stable normal bundle over  $M^2$. Define $\rho^{\odot}([(M^2, \eta_M, \Xi_M)]) \in V$ 
 by the formula:  $\rho^{\odot}(M^2, \eta_M, \Xi_M)=(M^2,\Xi_M)$, $[(M^2,\Xi_M)] \in V$.

\begin{lemma}\label{V}
The homomorphism $\rho^{\odot}: W_{\Q}^{\odot} \to V$ 
is decomposed as following:  
$$\rho^{\odot} = \rho \circ p: W_{\Q}^{\odot} \to W_{\Q} \to V,$$
where the homomorphism 
$\rho: W_{\Q} \to V$ is well-defined and is an epimorphism onto the index $2$ subgroup 
in $V$ of elements of the order $4$.
\end{lemma}

\subsubsection*{Proof of Lemma  $\ref{V}$}

Consider the standard transfer homomorphism with respect to the subgroup
$\Q \subset \aleph$, denote the transfer homomorphism by  $!: W \to W_{\Q}^{\odot}$. The following lemma is required.

\begin{lemma}\label{transf}
The image of the transfer homomorphism
$!: W \to W^{\odot}_{\Q}$ is inside the kernel  $Ker \rho^{\odot}$.
\end{lemma}

\subsubsection*{Proof of Lemma  $\ref{transf}$}
 A given arbitrary hyperquaternionic form $(M^2,\eta_M,\Xi_M)$, is represented by a connected surface. Take a geometrical stabilization of the surface $M^2$ by a connected sum with $2$ mirror copies of  Moebius bands, the generators of the bands are represented by
the element $T$ (we say that a band of the considered type is a $T$-band). Denote the result of the stabilization again by $(M^2,\eta_M,\Xi_M)$.
As the result, the surface $M^2$ is a connected sum of Moebius bands, which are represented by the elements
 $\j$, or by $\k$ (we say that a band of the considered type is a quaternion band). 
 
Take the decomposition of $M^2$ into a connected sum of Moebius bands with the only $T$-band and several quaternion bands.
 By the transfer homomorphism $M^2$ is covered by (a non-oriented) surface $\tilde M^2$. A $T$-band in the decomposition of $M^2$ is transformed into a cylinder on $\tilde M^2$, the generator $\tilde l \subset \tilde M^2$ of the cylinder is a closed loop on corresponds to the double covering over the generator $l \subset M^2$ of the $T$-band, the Hopf invariant $h(\tilde l) \in \Z/2$ of $\tilde l$ loop is trivial. 
 The each quaternion band on $M^2$ is covered by a pair of quaternion mirror-symmetric bands on $\tilde M^2$. This proves that the image   
 $(M^2,\eta_M,\Xi_M)^!$ in $V$ is trivial. Lemma  $\ref{transf}$ is proved. 
\[  \]

The last part of the proof of Lemma  
$\ref{V}$ is following. Let   $(M^2,\eta_M,\Xi_M)$ represents an arbitrary element in  $W_{\Q}$.
Consider the manifold 
 $P^3$ with boundary $\partial P^3 = M^2$, the manifold is equipped with a normal bundle structure
$(P^3,\zeta_P,\Psi_P)$, this structure determines a boundary of the  form  $(M^2,\eta_M,\Xi_M)$.
Denote by 
 $Q^2 \subset P^3$ a  surface, which is defined as a dual surface to  $\zeta_{\beta}$. Obviously, there exists a closed characteristic surface, because $\zeta_{\beta;P} \vert_{\partial P^3}$ is null-homotopic. Then
 $(Q^2, \zeta_P \vert_{Q}, \Psi_P \vert_Q)$ determines an element $x \in W$, the transfer $x^!$ belongs to  $Ker(\rho^{\otimes})$ by Lemma $\ref{transf}$. By the construction, $\rho^{\otimes}[(M^2,\eta_M,\Xi_M)]$ coincides with $x^!=(Q^2, \zeta_P \vert_{Q})^!$ in $V$. Lemma  $\ref{V}$ is proved.

\subsection*{Proof of Theorem  $\ref{th2}$}



Let us define a hyperquaternionic form
$(M^2,\eta_M, \Xi_M)$. Consider a pair of Moebious bands $(\mu_i,\partial) \subset M^2, i=1,2$, the generators of 
$\mu_1$, $\mu_2$ is represented by 
$\eta_M$ into the elements $TJ$, $T$ correspondingly. The connected sum $(\mu_1,\partial) \sharp (\mu_2,\partial)$ along
the common boundary $\partial \mu_1 = \partial \mu_2$ coincides to $M^2$. 
The characteristic mapping $\eta_M$ admits a reduction: $\eta_M: M^2 \to K(\D,1) \subset K(\aleph,1)$.

By the construction, $M^2$ contains a thin cylinder $C_J \subset M^2$, the (orientation preserved) loop $l_J \subset C_J$ which corresponds to the element $J \in \D \subset \aleph$ by $\eta_M$. The surface $M^2 \setminus C_J$ is diffeomorphic to the cylinder $C_{-J}$, this cylinder
is a non-oriented cycle between the two copies of $\partial C_J$. Denote the segment of the cylinder $C_{-J}$,
which is transversal to the central line of $C_{-J}$ by $l_T \subset C_{-T}$.
Extend the segment $l_T \subset C_{-T} \subset M^2$ by a closed loop on $M^2$ by a short path, which is transversal to $l_J$. This closed path is denoted by $l_{T} \subset M^2$. The closed path $l_{JT} \subset M^2$ are defined as
the central path in $M^2 \setminus l_{T}$. The paths $l_{T}$, $l_{JT}$ coincide with central lines the the Moebious bands
$\mu_1$, $\mu_2$ on $M^2$.

The  (orientation reversed) loop $l_T \subset M^2$ corresponds to the element $T \in \D \subset \aleph$ by $\eta_M$. The 
element $\eta_{M}(l_{T}^{-1} \circ l_J \circ l_{T} \circ l_J)$ 
is the trivial element in $\D \subset \aleph$, because $[T,J]=-1$. 
Informally speaking, the Klein bottle $M^2$ is the result of a non-oriented self-homology of $l_J$ by $l_T$.

Describe a regular cobordism of $2(M^2,\eta_M, \Xi_M)$ into a form $(L^2,\eta_L,\Xi_L)$, where $L^2$ is the 
Klein bottle, which is defined analogously to $M^2$. Take the orientation preserving loop  $l_{-1} \subset C_{-1} \subset L^2$,
which represents  the element $J^2=-1 \in \D \subset \aleph$ by $\eta_L$.
The loop $l_{-1}$  is the analog of the loop $l_J \subset M^2$.
Define the orientation reversed loop, which is analog of the loop $l_T \subset M^2$.
Denote the corresponding cycle of $l_{-1}$ by $S_1$, denote the corresponding cycle of $l_T$ by $S_2$.

\begin{lemma}\label{surgery1}
The form $2(M^2,\eta_M, \Xi_M)$ is equivalent to the form
 $(L^2,\eta_L,\Xi_L)$ (probably, up to an element of the order $2$ in $W$).
 \end{lemma}
 
Describe a regular cobordism of $2(L^2,\eta_L, \Xi_L)$ into a form $(K^2,\eta_K,\Xi_K)$, where $K^2$ is the 
Klein bottle, as in the case of $M^2$ and $L^2$. Denote the orientation preserved cycle  $R_2 \subset C \subset K^2$, which is the analog of the cycle $S_2 \subset C_{-1} \subset L^2$ and
which is represented into the trivial element in $\aleph$, by $\eta_K$. In this formula $C$ in a thin cylinder, which is the analog of the cylinder $C_{-1}$.
Denote the orientation reversed cycle $l_T \subset \mu_1 \subset K^2$ by $R_1$.

Recall that  an immersion
$f: K^2 \looparrowright \R^{10}$ with the prescribed isomorphism $\Xi_K: \nu_K \cong \eta_K^{\ast}(A) \oplus \eta_{\alpha;K}^{\ast}(\gamma) \oplus 3\eta_{\beta;K}^{\ast}(\gamma)$ of the normal bundle, where $A$ is the universal $4$-bundle over the subspace $K(\Z/2(-1) \oplus \Z/2(T),1) \subset K(\aleph,1)$, $\gamma$ is the universal line bundle,  is well-defined.
The element $\eta_K(R_1)$ is the trivial element in $\aleph$. Moreover, the mapping $\eta_K(R_1)$ has the target 
a point in $K(\aleph,1)$.

 The curve $f(R_1)$ is a framed curve in $\R^{10}$ and the stable Hopf invariant $h(R_1) \in \Z/2 = \{0,1\}$ is well-defined.

 \begin{lemma}\label{surgery2}
 The form   $2(L^2,\eta_L,\Xi_L)$ is equivalent, probably, up to an element of the order $2$ in $W$, to a form $(K^2, \eta_K, \Xi_K)$, where the oriented framed loop $R_1$ has the Hopf invariant $h(R_1) \ne 0, h(R_1) \in \Z/2$.
\end{lemma} 
 
\subsubsection*{Proof of Lemma $\ref{surgery1}$ and Lemma $\ref{surgery2}$}

Proofs of Lemmas are analogous. Let us prove  Lemma $\ref{surgery2}$.
The characteristic mapping $\eta_{L}$ takes the image in the subgroup $\Z/2(-1) \times \Z/2(T) \subset \aleph$,
where the generators of the factors are $\{-1,T\}$.

Define the normal bundle structure $\Xi_{L}$ as following. 
The normal bundle for  $(L^2,\eta_{L},\Xi_{L})$ is represented by a Whitney sum of 
4-bundle, 3-bundle and the trivial line bundle $A \oplus B \oplus \varepsilon$. 

The bundle $B$ is splitted into the Whitney sum of 3 isomorphic line bundles: $B=B_1 \oplus B_2 \oplus B_3$.
Each factor $B_j$, $j=1,2,3$ is the line bundle, which is skew along the cycle $R_2$ by means of the element $T$,
and is constant along the cycle $R_1$. The factors correspond to $\eta_{\beta;L}^{\ast}(\gamma)$.

The bundle $A$ is splitted into the  Whitney sum of 2 isomorphic copies of plane-bundles: $A=A_1 \oplus A_2$.
The plane bundle $A_{1}$ (and $A_{1}$) 
should be looked as a line complex bundle. The each line complex bundle is equipped with the Hermitian conjugation long the cycle $R_2$ by means of the point symmetry, given by multiplication on  $-1$  along the cycle $R_1$.  
 The  factors $A_{1},A_{2}$ are inside the 
4-dimensional block  $\eta_L^{\ast}(A)$ of $\nu_{L_2}$ with generators $\{\pm1,T\}$. 

The factor $\varepsilon$  corresponds to $\eta_{\alpha;L}^{\ast}(\gamma)$.


Denote two copies of $(L^2,\eta_{L},\Xi_{L})$ by $(L^2_1,\eta_{L_1},\Xi_{L_1})$, $(L^2_2,\eta_{L_2},\Xi_{L_2})$.
Define the following form $(L^2_2,\eta_{L_2}^{op},\Xi_{L_3}^{op})$, which represents an element in $W$.
The characteristic classes $\eta_{L_2}$, $\eta_{L_2}^{op}$
coincide, the normal bundle structure $\Xi_{L_2}^{op}$ is derived from  $\Xi_{L_2}$ by the reversing of the orientation
of the each factors $B=B_1 \oplus B_2 \oplus B_3$ and by 
the complex conjugation 
on the factors
$A_{1}, A_{2}$. In particular, the local orientations on the surfaces $(L^2_2,\Xi_{L_2})$ and $(L^2_2,\Xi_{L_2}^{op})$
with a prescribed normal bundle structure are opposite.

Let us prove that the forms $(L^2_2,\eta_{L_2},\Xi_{L_2})$,  $(L^2_2,\eta_{L_2},\Xi_{L_2}^{op})$
are equivalent in $W$. 
Take a self-homotopy of $\eta_{L_2}$ into itself such that the trace of a point
$pt \in L_2$ by this homotopy represents the generator $T \in \Z/2(-1) \times \Z/2(T) \subset \aleph$.
By this homotopy the framing $\Xi_{L_2}$ is transformed into a framing $\Xi_{L_2}^{op}$, where  $\Xi_{L_2}^{op}$ is the composition of
$\Xi_{L_2}$ with the reflection in the factors  $B_1, B_2, B_3, A_{1}, A_{2}$ as described above. The forms are equivalent.

Let us prove that the form $(L^2_1,\eta_{L_1},\Xi_{L_1}) \cup (L^2_2,\eta_{L_2},\Xi_{L_2}^{op})$ 
is regular cobordant to the form $(K^2,\eta_K,\Xi_K)$, probably, up to a form $(P^2,\eta_P,\Xi_P)$ with the characteristic class $\eta_P$ takes the image in the central subgroup $\Z/2(-1) \subset \aleph$.



Take the restriction of $\Xi_1$ and of $\Xi_2^{op}$ over the cycle $R_1 \subset L_1^2$ and the cycle $-R'_1 \subset L^2_3$ correspondingly
(in this formula $-R'_1$ is the cycle on $L^2_2$ which corresponds to the  cycle   $R_1$ with the opposite orientation, using the diffeomorphism $L^2_1 \cong L^2_2$).
The restrictions  $\Xi_1 \vert_{R_1}$, $\Xi_2^{tw} \vert_{-R_1}$ are 4-dimensional $(-1,T)$-framings, which are stabilized in the codimension $4$ by corresponding framings on $B \oplus \varepsilon$ ($B_i \vert_{R_1}, i=1,2,3$ is the trivial bundle, the trivialization $\Xi_2$  over $B\vert_{R_1}$ is opposite to the trivialization  $\Xi_2^{tw} \vert_{R_1}$, the trivialization $\Xi_2$  over $A\vert_{R_1}$ is conjugate to the trivialization  $\Xi_2^{tw} \vert_{R_1}$).

Denote a $(-1,T)$-framing $\Xi_2^{tw}$ over $(L_2,\eta_{L_2})$ as following. 
Denote the line subbundles $\lambda_1 \subset A_1$, $\lambda_2 \subset A_2$, which correspond to the imaginary axis of the complex line bundles.
The line bundle $\lambda$ over $L_2^2$ is well-defined, and this bundle is skew over the cycle $R_2 \subset L_2^2$,
which corresponds to the element $T$.  
Take the rotation trough the angle $\pi$ inside the 4-bundle $B \oplus \lambda$ over $L_2$. 
As the result we get the new $(-1,T)$-framing over $L_2^2$, denoted by $\Xi_2^{tw}$.
The framing $\Xi_2^{tw}$ coincides to the framing $\Xi_1$ everywhere, except the line bundle $\lambda_2 \subset \nu_{L_1}$, on this factor the framing $\Xi_2^{tw}$ is given by the reflection of $\Xi_1$. 
The framings $\Xi_2^{tw}$ is equivalent to the framing $\Xi_2^{op}$, and is equivalent to the framing $\Xi_1$.

Assume without loss of a generality that the restriction of the framing $\Xi_1 \vert_{R_1}$ to the subbundle $B \oplus \varepsilon$ over the cycle $R_1$ is parallel to the coordinate axis
$e_7,e_8,e_9,e_{10}$. Assume the framing $\Xi_1 \vert_{R_1}$ on the factors $\lambda_1, \lambda_2$ is parallel to the 
vectors $e_4,e_6$ correspondingly. Assume the framing $\Xi_1 \vert_{R_1}$ on the factors $A_1, A_2$ is parallel to the 
vectors $(e_3,e_4),(e_5,e_6)$ correspondingly. Then the skew-framing $\Xi_2^{tw} \vert_{R_1}$ coincides to
the $\Xi_1$ along each directions, but the direction of the coordinate vector $e_6$, where  $\Xi_1$, $\Xi_2^{tw}$ are opposite.

The $-1$-structure of skew framings $\Xi_1 \vert_{R_1}$, $\Xi_2^{tw} \vert_{-R'_1}$ are distinguished only
inside the factor $A_2$ of the normal bundle of $L_1^2 \cong L_2^2$, by a
full rotation trough the angle $2\pi$. 

Take the regular cobordism transformation of the form $(L^2_1,\eta_{L_1},\Xi_{L_1}) \cup (L^2_2,\eta_{L_2},\Xi_{L_2}^{tw})$
by a surgery, with a support in small neighborhoods of a corresponding pair of points on $R_1$, $-R'_1$. As the result we get the form $(L^2_4,\eta_{L_4},\Xi_{L_4})$.
The image of $\eta_{L_4}$ is in the space $K(\Z/2(-1) \times \Z/2(T),1)$.
The cycle $R_1 \cup -R'_1$ is transformed into a cycle $R_3 \subset L^2_3$. The image of the characteristic class $\eta_{L_3}(R_3)$
is null-homotopic in the target space  $K(\Z/2(-1) \times \Z/2(T),1)$. The stable Hopf invariant $h(R_3)$ of the framed curve $R_3$ is non-trivial. 

The form $(L^2_3,\eta_{L_3},\Xi_{L_3}) \cup (K^2, \eta_K, \Xi_K)$ is cobordant to a form  
$(L^2_4,\eta_{L_4},\Xi_{L_4})$, where the image of the mapping $\eta_{L_4}$ is inside the space
$K(\Z/2(-1),1)$. The form $(L^2_4,\eta_{L_4},\Xi_{L_4})$ is trivial, or, is of the order $2$ in $W$.
Theorem 
$\ref{th2}$ is proved.

\section*{Conclusion} 
V.I.Arnol'd formulated the  problem [\cite{Arn}, Problem 
1984-12]: "To transform asymptotic ergodic definition of the Hopf invariant
of divergence-free vector fields to the theory of S.P.Novikov, which 
generalize the Whitehead product of homotopy groups of spheres"'.

Algebraic commutators, which are used to define the higher invariants   of classical links, are particular 
Whitehead products in homotopy groups of spheres. $M$-invariant  is a special generalized Whitehead product, which admits asymptotic and ergodic property. To keep additional symmetry of magnetic fields we have to 
apply the $M$-invariant for links in various homogeneous manifolds, which are rational homology spheres. For classical links $M$-invariant is associated with the Arf-invariant in the stable homotopy group $\Pi_2$. Hypothetic modifications of $M$-invariant for links in $S^3/\Q$ are associated with Arf-Brown invariant in the stable homotopy group $\Pi_3$, and with hyperquaternionic Arf-invariant in the stable homotopy group $\Pi_7$. The constructions give a solution (in part) of the Arnol'd Problem.

\end{document}